\documentclass{amsart} 
\usepackage{amsmath}
\usepackage{amssymb}
\usepackage{amsthm}

\newtheorem{thm}{Theorem}[section]
\newtheorem{lem}[thm]{Lemma}
\newtheorem{prop}[thm]{Proposition}

\newtheorem{defn}[thm]{Definition}
\numberwithin{equation}{section}

\newcommand{\R}{{\mathbb R}}

\newcommand{\Z}{{\mathbb Z}}

\newcommand{\N}{{\mathbb N}}

\begin{document}

\title{ Completions of altered topological subgroups of $\R^n$}
\author{Jon W. Short}
\subjclass[2000]{Primary 22A05; Secondary 54G15, 54E35}
\address{Department of Mathematics and
Statistics, Sam  Houston State University,\linebreak Huntsville,
Texas 77341}
\email{jon@shsu.edu} 
\urladdr{http://www.shsu.edu/$\sim$\;mth\_jws}
\keywords{topological group, totally disconnected,
local isometry, not locally compact, convergent sequence of real numbers}


\begin{abstract}
We prove that a large class of metrizable group topologies for subgroups of $\R^n$ and the completions of the subgroups are locally isometric to, respectively, metrizable group topologies for $\Z$ and their completions, first studied by Nienhuys.  This will  prove, in particular, that all the complete groups in question are one dimensional, locally totally disconnected, and not locally compact.  The metrizable topologies on the subgroups of $\R^n$ are formed by specifying a sequence in $\R^n$ and the rate at which it must converge to the identity.
\end{abstract}

\maketitle

\section{Introduction}
In this paper, we examine the local structure of complete topological groups obtained by manipulating the topologies of certain subgroups of $\R^n$.  These groups are constructed using pairs of sequences $(\{v_j\},\{p_j\})$ where $\{v_j\}$ is a sequence in $\R^n,$ $\{p_j\}$ is a sequence of positive real numbers, and  a primary property of the pair is that $\|v_j\|\rightarrow \infty$ and $p_j\rightarrow 0$ in the standard topologies.  A specific example of such a pair is $(\{j!+\sqrt{2}\},\{1/j\}).$  The groups in question are the  subgroups of $\R^n$ generated by $\{v_j\}$ (which will be denoted $G_v$) and the group topologies are obtained by  ``forcing'' $\{v_j\}$ to converge to zero at least as fast as $\{p_j\}$ converges to zero in the standard topology.    Since the groups $G_v$ are generally not complete, we then form the completions  and study the properties of these groups.  (We will prove that the topologies of these groups will be either weaker or not comparable to the standard topologies.)  The details of this construction are given in \cite{lclisom} and first studied in \cite{cs:ds}.   Our primary goal in this paper is to show that  many of the local properties of these groups can be understood.  

Other authors have certainly used pairs of sequences to obtain interesting results on the structure of complete topological groups.  In particular, J.W. Nienhuys, in a group of papers \cite{nh2, nh1, nh3}, studies complete monothetic groups obtained using pairs of sequences $(\{v_j\},\{p_j\})$ where $\{v_j\}$ is a sequence of natural numbers and $\{p_j\}$ is a sequence of real numbers.  With certain restrictions placed on the pair of  sequences, he is able to give information regarding the  topological structure of the complete monothetic groups, as well as examples of topological groups exhibiting unusual properties.  We show in this paper that we can relate via a local isometry or local homeomorphism the groups considered by Nienhuys and those constructed by this author and Stevens in \cite{lclisom}.  The process of obtaining this local isometry will include creating a ``special'' sequence that bridges the gap between the sequences used by Nienhuys and those used in \cite{lclisom}.    Thus we are able to exploit Nienhuys's analysis of the completions of unusual topologies on the integers to obtain a deeper understanding of the structure of the completions of unusual topologies on subgroups of $\R^n.$

This paper is divided into several sections with notation and background material given in \S2 and 3 and the main results  given in \S 4 and 5.  Examples are provided throughout the paper.

\section{Notation and terminology} 
We will use the
notation presented in this section throughout the paper. 
Other notation will be introduced as needed.

We denote the topological space  $ G$ with topology 
$\mathcal{T}$ as $(G, \mathcal{T})$.  If $ d$ is a metric
on a topological space $ G,$ or a norm that induces a
metric, we will blur the distinction between $d$ and
the topology it induces and denote the space with the
metric topology by
$(G, d)$.  

All topological groups in this paper are abelian (we
will denote the identity element by $0$), metrizable, and separable. Since the
groups  are metrizable, a natural completion to consider is the group's completion as a metric space.  Since the groups are abelian, the group operations can be extended in a continuous way to that completion.  We will denote that
completion of the topological group $ (G,
\mathcal{T})$ by $\mathcal{C}(G,
\mathcal{T})$.  Also, since the groups are separable metric spaces, all dimensions functions will coincide for these groups.

$ \N, \Z,$ $\R$, and $\R^{+}$ will denote, respectively, the natural
numbers, the integers, the real numbers, and the positive real numbers. The standard
topology on $ \R^j$ will be denoted $ \tau^j$ and if
$j=1$ as $\tau$.  If
$ x\in {\R}^{m},$ for $ m\in\N,$ then $
\|x\|$ will be the standard norm of $ x$.  If $ m=1,$ $
\|x\|$ will be abbreviated $ |x|$.  If $ x\in\R,$ then $
\left\lceil x \right\rceil$ will denote the greatest
integer less than or equal to $ x$.

All sequences in a topological space will be indexed on $
\N$. Let
$ \{x_j\}$ and $\{y_j\} $ be  sequences in $ {\R}^{m}$ and
$ {\R}^{n},$ respectively.  We will often have cause to
construct a sequence in $ {\R}^{m+n}$ from $ \{x_j\}$ and
$\{y_j\}$.  To see how this is accomplished, let $
x_j=(x_{j1}, x_{j2},\ldots, x_{jm})$ and  $y_j=(y_{j1},
y_{j2},\ldots,y_{jn})$.  Then $$
\{(x_{j1}, x_{j2},\ldots, x_{jm},y_{j1}, y_{j2}, \ldots,y_{jn}): j\in\N\}$$ will be a sequence in $
{\R}^{m+n}$.  We will denote this sequence by
$\{(x_j,y_j)\}$.

\section{Background Information} 

In \cite{lclisom} we show that large classes of metrizable groups, obtained by adjusting the topology on the standard Lie groups $\R^n$, are locally isometric.  This allows us to understand the local nature of many topological groups by understanding the local structure of one.  To make this idea precise, we first need a definition.

\begin{defn}\mbox{}

\begin{enumerate}
\renewcommand{\labelenumi}{(\roman{enumi})}
\item Let $
\{p_j\}$ be a nonincreasing sequence of numbers in $\R^{+}$, which converges to zero in the standard
topology on $
\R$, and let $ \{v_j\}$ be a sequence in $ \R^n-\{0\}$ such that $ \{\|v_j\|\}$ is
nondecreasing and the sequence $
\{p_{j+1}\|v_{j+1}\|/\|v_j\|\}$ has a positive lower
bound.  The pair $(\{v_j\},\{p_j\})$ will be called a sequential norming pair for $\R^n$ and abbreviated (SNP).  
\item If $(\{v_j\},\{p_j\})$ is a (SNP) for $\R^n$ and $\{q_j\}$ is a sequence in $\R^m$, then the pair $(\{(v_j,q_j)\},\{p_j\})$ will be called an extended norming pair for $\R^{n+m}$ and abbreviated (ENP).
\end{enumerate}
\end{defn}

An example of a (SNP) for $\R$ is $(\{j!+\sqrt{2}\},\{1/j\})$.  From this we can easily obtain an (ENP) for $\R^n$.  As a specific example, $(\{(j!+\sqrt{2},j)\},\{1/j\})$ is an (ENP) for $\R^2$.  We note that not every (SNP) for $\R^n, \; n>1$ is obtainable as a (ENP).  To see this, consider the (SNP) for $\R^2$ given by $(\{1/j, j!\},\{1/j\})$.  It is not a (ENP) since it would have to be extended from the pair $(\{1/j\},\{1/j\}),$ which is not a (SNP).  Likewise, not every (ENP) is a (SNP).  In fact, (ENP)s differ from (SNP)s in a significant way.  If  $(\{v_j\},\{p_j\})$ is a (SNP), then the condition that $\{p_{j+1}\|v_{j+1}\|/\|v_j\|\}$ has a positive lower
bound implies that the sequence $\{\|v_{j+1}\|/\|v_j\|\}$ diverges to infinity.  If $\{w_j\}=\{(v_j,q_j)\}$ is part of a corresponding (ENP), however, then the sequence $\{\|w_{j+1}\|/\|w_j\|\}$ does not necessarily diverge to infinity (see \cite[pg. 50]{lclisom} for a specific example), although it does have to be unbounded.  We note this difference since the quotients $\{\|v_{j+1}\|/\|v_j\|\}$ and $\{\|w_{j+1}\|/\|w_j\|\}$ play an important role in the following exposition and quotients of this type are also important to those interested in the study of precompact group topologies (see \cite{precompact} for example).

To get our ``unusual'' subgroups that we are going to study, we start with a (SNP) or (ENP) $(\{v_j\},\{p_j\})$ and force $\{v_j\}$ to converge to the identity in the subgroup generated by $\{v_j\}$ using the following norm.     
 
 We assume that all sums are finite unless otherwise specified.
\begin{thm}\label{normthm}\cite[pg 54]{lclisom}
 Let $(\{v_j\}, \{p_j\})$ be a
(SNP) on $\R^{n}$, and let $G_v$ be the subgroup of
$\R^{n}$ generated by the sequence $\{v_j\}$. Then 
$\sigma:G_v\rightarrow \R$ defined by
$$\sigma(x)=\inf \left\{ \sum{|c_j|p_j}:x=\sum{c_j v_j,\; c_j\in\Z}
\right\}$$ induces a metrizable group topology on $G_v$ in
which
$\{v_j\}$ converges to zero. If $\{q_j\}$  is any sequence
in $\R^m$, $\{w_j\}=\{(v_j, q_j)\}$, and $G_w$ is the
subgroup of $\R^{n+m}$ generated by $\{w_j\}$, then 
$$\sigma'(x)=\inf \left\{ \sum{|c_j|p_j}:x=\sum{c_j w_j,\; c_j\in\Z}
\right\}$$ induces a metrizable group topology on $G_w$ in
which
$\{w_j\}$ converges to zero.
\end{thm}

Although $\sigma$ and $\sigma'$ look similar, they are in general, operating on different spaces.   For example, if we consider the norming pair $(\{j!+\sqrt{2}\},\{1/j\})$ and the (ENP) $(\{(j!+\sqrt{2},j!)\},\{1/j\})$, then $\sigma$ is a function with domain a subgroup of $\R$ while $\sigma'$ is a function with domain a subgroup of $\R^2$.  At this point we can see that (ENP)s are useful for easily obtaining altered topological  subgroups of $\R^n.$  All we need is a (SNP) for $\R$ and we can get (ENP)s, and their corresponding altered subgroups, in $\R^n$ for any $n>1.$

We are now in a position to state a main result from \cite{lclisom}.

\begin{prop}\label{thm:isom}\cite[Proposition 14]{lclisom}
Let $ (\{v_{j}\},\{p_{j}\})$
and $ (\{n_{j}\},\{p_{j}\})$ be two (SNP)s on $\R^{n}$
such that
\begin{equation*}
p_j\left\lceil\frac{\|v_{j+1}\|}{\|v_j\|}-1\right\rceil,
p_j\left\lceil\frac{\|n_{j+1}\|}{\|n_j\|}-1\right\rceil
\geq 1
\end{equation*} for all but finitely many $j,$ and let $
(\{v_{j}\}',\{p_{j}\})$ and $ (\{n_{j}\}',\{p_{j}\})$ be
two extended norming pairs  corresponding to,
respectively, $ (\{v_{j}\},\{p_{j}\})$
and $ (\{n_{j}\},\{p_{j}\})$.  Let
$(G_v,
\sigma_v)$ be the subgroup of $\R^{n}$  generated by
$\{v_j\}$ with topology induced by the norm $\sigma$ given
in
\ref{normthm}.  $(G_n, \sigma_n)$, $(G_{v'},
\sigma_{v'})$, and
$(G_{n'}, \sigma_{n'})$ are defined similarly.  Then
\begin{enumerate}
\renewcommand{\labelenumi}{(\roman{enumi})}
\item $(G_v, \sigma_v)$ is locally isometric to $(G_n,
\sigma_n)$.
\item $(G_{v'}, \sigma_{v'})$ is locally isometric to
$(G_{n'}, \sigma_{n'})$.
\item $\mathcal{C}(G_v, \sigma_v)$ is locally isometric to
$\mathcal{C}(G_n, \sigma_n)$.
\item $\mathcal{C}(G_{v'}, \sigma_{v'})$ is locally
isometric to
$\mathcal{C}(G_{n'}, \sigma_{n'})$.
\end{enumerate}
\end{prop}

Above we saw that $(\{(j!+\sqrt{2},j)\},\{1/j\})$ is an (ENP) for $\R^2$.  We can easily get a different (ENP) for $\R^2$ by replacing the second coordinate of the sequence.  For example, $(\{(j!+\sqrt{2},1/j)\},\{1/j\})$ is another (ENP) for $\R^2$.  Yet another (ENP) for $\R^2$ is $(\{(2^{2^j},1)\},\{1/j\}).$  These sequences exhibit very different behaviors in the standard sense and generate different topological groups, but \ref{thm:isom} shows that the local topological structure of the topological groups is the same.

With this theorem at our disposal, we ask if the completion of one of these topologies is understandable.  If so, we can then understand the local structure of many.  Results by Nienhuys will help us in this case.

Nienhuys also uses pairs of sequences $(\{v_j\},\{p_j\})$ to generate group topologies, but in this case $\{v_j\}$ is a sequence of natural numbers and the pair satisfies certain other properties.

\begin{defn}
The pair of sequences $(\{v_j\},\{p_j\})$ will be called a Nienhuys norming pair, and abbreviated (NNP), if $v_1=1,\; v_j\in\N,\; v_j|v_{j+1},\; p_j\in\R^{+}$ and $p_{j+1}v_j\leq p_j v_{j+1}.$

\end{defn}
The condition that $v_j|v_{j+1}$ is very important in Nienhuys's arguments and allows him to obtain some very nice results concerning the completions of  group topologies on $\Z$.  We note that none of the previously mentioned (SNP)s have the property that $v_j|v_{j+1}$.    Some are not even sequences of integers.  So we see that there is a clear distinction between (SNP)s and (NNP)s, even when the (SNP) generates a subgroup of $\R.$  (Recall that (SNP)s generate subgroups of $\R^n$ for $n\geq 1.$)  On the other hand, there are also (NNP)s that are not (SNP)s.  To see this, consider the pair $(\{2^{j-1}\},\{1/j\}).$  It is clearly a (NNP), but it is not a (SNP) since
$$\left\{\frac{1}{j+1}\cdot\frac{2^j}{2^{j-1}}\right\}=\left\{\frac{2}{j+2}\right\}$$ does not have a positive lower bound.

The next theorem summarizes Nienhuys's results that we will use to study our topological groups.  We need, however, one definition before proceeding.

\begin{defn}
A topological space is totally disconnected if its components  are points.

\end{defn}

\begin{thm}\label{thm:nnp}\cite[Lemmas 6,7, Theorems 24, 67]{nh1}
Let $(\{v_j\},\{p_j\})$ be a (NNP) and define $k_j=\frac{v_{j+1}}{v_j}$.  
\begin{enumerate}
\renewcommand{\labelenumi}{(\roman{enumi})}
\item There is a norm $\nu$ on $\Z$ such that $\nu(v_j)=p_j$, and the largest such norm is given by $\nu(x)=\inf\{\sum{|c_j|p_j}\;:\; x=\sum{c_j v_j},\; c_j\in\Z\}$.
\item $\mathcal{C}(\Z,\nu)$ is compact if and only if $\sum_{j=1}^{\infty}k_jp_j<\infty,$ is Hausdorff if and only if $p_j>0\; \forall j$, and is nondiscrete if and only if $\lim\inf p_j=0$.
\item $\mathcal{C}(\Z,\nu)$ is zero dimensional if and only if it is locally compact, and is one dimensional if and only if it is not locally compact.  In any case, $\mathcal{C}(\Z,\nu)$ is totally disconnected.

\end{enumerate}

\end{thm}

We stated the results above in the same manner that Nienhuys stated them.  We do note, however, that in this situation, locally compact and compact are equivalent terms \cite[Theorem 9.1]{hr}. Be clear that even though the  ``largest'' norm mentioned in one above, and the  norm described in \ref{normthm} look very similar, they are (most likely) defined on entirely different spaces.

\section{Moving (SNP)s into the framework of (NNP)s}

In this section we will begin with a (SNP)  and show how it can be associated with a corresponding (NNP).  This correspondence will allow us to examine the local nature of the topological groups and their completions which are determined by the (SNP).  Note that if $(\{v_j\},\{p_j\})$ is a (SNP), then $\|v_j\|$ will eventually be greater than one.  Since we are ultimately interested in the local topological structure of groups determined by (SNP)s, and since the same local topological structure will be generated if we ignore finitely many terms of the sequence,  we will assume that $\|v_j\|\geq 1\; \forall j\in\N.$

\begin{lem}\label{lem}
Let $(\{v_j\},\{p_j\})$ be a (SNP)  for $\R^n$.  Define $x_1=1$, and for $j\geq 2,$ define
$$x_j=\lceil \|v_1\| \rceil \lceil \|v_2\| \rceil\cdots \lceil \|v_j\| \rceil.$$
\begin{enumerate}
\renewcommand{\labelenumi}{(\roman{enumi})}
\item The pair $(\{x_j\},\{p_j\})$ is a (SNP)  for $\R$ and a (NNP) for $\Z$.
\item If $(\{v_j\},\{p_j\})$ satisfies the hypotheses of \ref{thm:isom}, then so does $(\{x_j\},\{p_j\}).$

\end{enumerate}

\end{lem} 
\begin{proof}
Recall that $\|v_j\|\geq 1\; \forall j\in\N.$  Since $\{\|v_j\|\}$ is nondecreasing, $\{\|x_j\|\}$ is nondecreasing.  All that remains to prove (i) is to show that the sequence $\left\{p_{j+1}\frac{\|x_{j+1}\|}{\|x_j\|}\right\}$ has a positive lower bound.  Since $\|v_j\|\geq 1\; \forall j\in\N$, we know that $\|v_{j+1}\|/\lceil \|v_{j+1}\| \rceil \leq 2$.  Since ${\|v_j\|}\rightarrow \infty$, there is a $k\in\N$ such that for all $j\geq k,$ $\|v_j\| \geq 2.$  Thus, for $j\geq k$, 
$\lceil \|v_{j+1}\| \rceil \geq \frac{\|v_{j+1}\|}{\|v_j\|}.$  Now
$$p_{j+1}\frac{\|x_{j+1}\|}{\|x_j\|}=p_{j+1}\lceil\|v_{j+1}\|\rceil\geq p_{j+1}\frac{\|v_{j+1}\|}{\|v_j\|}.$$
Since $(\{v_j\},\{p_j\})$ is a (SNP), there is an $\alpha\in\R,\; \alpha>0$
 such that $p_{j+1}\frac{\|v_{j+1}\|}{\|v_j\|}\geq \alpha$.  Thus, for $j\geq k$, $p_{j+1}\frac{\|x_{j+1}\|}{\|x_j\|}\geq\alpha$.  If $j<k$, then $p_{j+1}\lceil\|v_{j+1}\|\rceil>0$ since $\lceil\|v_{j+1}\|\rceil\geq 1$ and $p_j>0 \;\forall j\in\N$.  Therefore $\left\{p_{j+1}\frac{\|x_{j+1}\|}{\|x_j\|}\right\}$ has a positive lower bound and $(\{x_j\},\{p_j\})$ is a (SNP).  Now we need to show that $(\{x_j\},\{p_j\})$ is also a (NNP).
 Using that $\|v_{j+1}\|\geq 1$, we see that
 $$\frac{p_{j+1}}{p_j}\leq 1\leq \lceil\|v_{j+1}\|\rceil=\frac{\|x_{j+1}\|}{\|x_j\|}.$$
 This along with the fact that $x_1=1, x_j\in\N,$ and $x_j|x_{j+1}$ proves that $(\{x_j\},\{p_j\})$ is a (NNP).
 
 To prove (ii), all we need to show is that 
$$p_j\left\lceil \frac{x_{j+1}}{x_j}-1\right\rceil \geq 1$$ for all but finitely many $j\in\N.$
Since $(\{v_j\},\{p_j\})$ satisfies the hypotheses of  \ref{thm:isom}, $p_j\left\lceil \frac{\|v_{j+1}\|}{\|v_j\|}-1\right\rceil \geq 1$ for all but finitely many $j\in\N$. Now
$p_j\left\lceil \frac{x_{j+1}}{x_j}-1\right\rceil = p_j\left(\lceil\|v_{j+1}\|\rceil - 1\right)$.  Since $(\{v_j\},\{p_j\})$ is  a (SNP), there is a $k\in\N$ such that for $j\geq k$ $\lceil \|v_{j+1}\| \rceil \geq \frac{\|v_{j+1}\|}{\|v_j\|}$, and thus for $j\geq k$, $\lceil \|v_{j+1}\| \rceil -1\geq \frac{\|v_{j+1}\|}{\|v_j\|}-1.$
Therefore \begin{eqnarray*}
p_j\left(\lceil\|v_{j+1}\|\rceil - 1\right) & \geq & p_j\left(\frac{\|v_{j+1}\|}{\|v_j\|}-1\right) \\
& \geq & p_j\left\lceil\frac{\|v_{j+1}\|}{\|v_j\|}-1\right\rceil\\
& \geq & 1
\end{eqnarray*}
for all but finitely many $j\in\N$.  Therefore $(\{x_j\},\{p_j\})$ satisfies the hypotheses of \ref{thm:isom}.
\end{proof}

To illustrate the previous theorem, consider the   (SNP) for $\R$ given by $(\{j!+\sqrt{2}\},\{1/j\})$.  The corresponding (NNP) is $(\{x_j\},\{1/j\})$ where 
$x_1=1, x_2=1\cdot (2!+1), x_3=1\cdot (2!+1)\cdot (3!+1), \ldots, x_j=1\cdot (2!+1)\cdot (3!+1)\cdots (j!+1).$  In this case, the (SNP) $(\{j!+\sqrt{2}\},\{1/j\})$ satisfies the hypotheses of \ref{thm:isom}, and thus we know by the previous lemma that $(\{x_j\},\{1/j\})$ does as well.

For the remainder of this section we will assume that all (SNP)s $(\{v_j\},\{p_j\})$ have the property that $\{v_j\}$ is a sequence in $\R$.  We will deal with the case of $\{v_j\}$ being a sequence in $\R^n$ for $n>1$, as well as (ENP)s,  in the next section. 

\begin{prop}
Suppose that $(\{v_j\},\{p_j\})$ is a (SNP) satisfying the hypotheses of \ref{thm:isom} and that $\{x_j\}$ is the sequence defined in \ref{lem}.  If $(G_v,\sigma_v)$ is the topological subgroup generated by $\{v_j\}$ and $(\Z,\sigma_x)$ is the  topological subgroup generated by $\{x_j\}$, then $(G_v,\sigma_v)$ is locally isometric to $(\Z,\sigma_x),$ and $\mathcal{C}(G_v,\sigma_v)$ is locally isometric to $\mathcal{C}(\Z,\sigma_x).$

\end{prop}
\begin{proof}
By \ref{lem}  $(\{x_j\},\{p_j\})$ is a (SNP) for $\R$ satisfying the hypotheses of \ref{thm:isom}. Thus $(\{v_j\},\{p_j\})$ and $(\{x_j\},\{p_j\})$ satisfy the hypotheses of \ref{thm:isom}, which gives the desired result.
\end{proof}

\begin{thm}\label{mres}
Suppose that $(\{v_j\},\{p_j\})$ is a (SNP) satisfying the hypotheses of \ref{thm:isom}.  If $(G_v,\sigma_v)$ is the topological subgroup generated by $\{v_j\}$, then $\mathcal{C}(G_v,\sigma_v)$ is one dimensional, locally totally disconnected, and not locally compact.
\end{thm}
\begin{proof}
$(\{x_j\},\{p_j\})$ is both a (SNP) and (NNP).  Since it is an (SNP), 
$$p_{j+1}\frac{x_{j+1}}{x_j}\geq\alpha >0$$
 for all $j\in\N$.  Thus $p_j\frac{x_{j+1}}{x_j}\geq \frac{p_j}{p_{j+1}}\alpha$.  Since $\{p_j\}$ is a decreasing sequence,  $\frac{p_j}{p_{j+1}}\alpha\geq \alpha >0.$  Therefore $\sum_{j=1}^{\infty}p_j\frac{x_{j+1}}{x_j}$ is unbounded.  Since $(\{x_j\},\{p_j\})$ is also a (NNP), we apply \ref{thm:nnp} to conclude that  $\mathcal{C}(\Z,\sigma_x)$ is not locally compact, is one dimensional, and is totally disconnected.  Since $\mathcal{C}(G_v,\sigma_v)$ is locally isometric to $\mathcal{C}(\Z,\sigma_x)$, we have the statement.
\end{proof}

Note that we are not claiming that $\mathcal{C}(G_v,\sigma_v)$ is totally disconnected.  We only know, due to the local nature of \ref{thm:isom}, that it has a neighborhood of the identity which is totally disconnected.  This is enough, however, to conclude that $\mathcal{C}(G_v,\sigma_v)$ is totally pathwise  disconnected.\footnote{Thanks to professor Christine Stevens for suggesting this result.}  

\begin{prop}
$\mathcal{C}(G_v,\sigma_v)$ contains no nontrivial paths.

\end{prop}
\begin{proof}Suppose the contrary.  Then there would be a nontrivial path $f:[0,1]\rightarrow \mathcal{C}(G_v,\sigma_v)$ with $f(0)=0$.   Let $c=\inf\{t\in [0,1]: f(t)\neq 0\}$.  Since the path is nontrivial, $c<1$.  We have two cases to consider.  First suppose that $c>0$.  Let $B$ be a totally disconnected neighborhood of the identity.  Then $f(c)+B$ is a totally disconnected neighborhood of $f(c)$.  Now $f^{-1}(f(c)+B)$ is open in $[0,1]$ and contains $c$.  Hence $f^{-1}(f(c)+B)$ contains an interval of the form $[a,b]$ where $a<c<b$.  Now $f([a,b])$ would be a nontrivial connected subset of $f(c)+B$ since $f(a)=0$ and $f(b)\neq 0$.  This is clearly a contradiction.   Now suppose that $c=0.$   Then since $f^{-1}(B)$ is open in $[0,1]$ and contains $c=0$, it contains an interval of the form  $[0,b]$ where $f(b)\neq 0.$  Thus $f([0,b])$ is a nontrivial connected subset of $B,$ and this is a contradiction.  Hence $\mathcal{C}(G_v,\sigma_v)$ contains no nontrivial paths.
\end{proof}

The question of whether $\mathcal{C}(G_v,\sigma_v)$ is totally disconnected remains open.  We note, however, that if $\mathcal{C}(G_v,\sigma_v)$ is not totally disconnected, then we can prove false a conjecture by A. Gleason.  In particular, Gleason has conjectured that a separable, metrizable, connected, topological group of finite dimension has a dense arc component of the identity, if the group is also complete \cite[pg. 296]{cs:wag}.  If $\mathcal{C}(G_v,\sigma_v)$ is not totally disconnected, then the completion of the connected component of the identity in $\mathcal{C}(G_v,\sigma_v)$  provides a counter example to Gleason's conjecture, since it would be complete, connected, separable, metrizable, and finite dimensional but, as the previous proposition shows, contain no nontrivial paths.

We are now in a position to describe how the groups that we are studying relate to the groups with the standard topologies.  We note that $\{v_j\}$ will generate either a dense or discrete subgroup of $(\R,\tau)$.  
If the subgroup  generated by $\{v_j\}$ is discrete, then $(G_v,\sigma_v)$ is strictly weaker than the standard topology on $G_v$ \cite[Proposition 12]{lclisom}.  On the other hand, if the subgroup generated by $\{v_j\}$ is dense, then $(G_v,\sigma_v)$ and $(G_v,\tau)$ are not comparable. 

\begin{prop}
If $\{v_j\}$ generates a dense subgroup of $(\R,\tau),$ then $(G_v,\sigma_v)$ and $(G_v,\tau)$ are not comparable.
\end{prop}
\begin{proof}
If  $(G_v,\tau)$ and $(G_v,\sigma_v)$ are comparable, then $(G_v,\sigma_v)$ is necessarily weaker than $(G_v,\tau)$, since it is has a $\sigma_v$-convergent sequence that is not $\tau$-convergent.  Given this, $\mathcal{C}(G_v,\tau)=(\R,\tau)$ would map into $\mathcal{C}(G_v,\sigma_v).$  This is a contradiction since the previous proposition  shows that $\mathcal{C}(G_v,\sigma_v)$ cannot contain nontrivial paths.
\end{proof}

\section{Subgroups Generated by (SNP)s in $\R$ are Fundamental}
Recall that (SNP)s and (ENP)s allow us to alter the topological group structure of certain subgroups of $\R^n$.  Given the very different nature of (SNP)s and (ENP)s, it is reasonable to assume that they might generate groups with different local properties.  On the other hand, (ENP)s are formed from (SNP)s, and thus it is also reasonable to assume that the groups generated are ``associated.''  In this section, we show that the second assumption is correct.  We do this by  showing that understanding the local structure of a group generated by a (SNP) in $\R$ is fundamental to understanding the local properties of subgroups generated by (SNP)s in $\R^n.$  In particular, we will prove that if $(\{v_j\},\{p_j\})$ is a (SNP) for $\R^n,\; n>1,$ then the topological group generated by $\{v_j\}$ is locally homeomorphic to a  group generated by a (SNP) in $\R.$  A similar statement is true for an (ENP) in $\R^{n+m}$, i.e., every (ENP) in $\R^{n+m}$ generates a group that is locally homeomorphic to a group generated by an (SNP) in $\R^n$.  With this we can then use results from the last section to understand the local topological structure of these groups.

At this point, we caution the reader  not to assume that even though we will prove that the groups formed by (SNP)s in $\R^n$ and (ENP)s are locally homeomorphic to groups formed by (SNP)s in $\R,$ that the extended groups are somehow uninteresting.  We can easily see in fact that the groups formed from (ENP)s are not all topologically isomorphic, and even if the groups are the same, they do not necessarily have the same topology.  For example, consider the (SNP) given by $(\{j!\},\{1/j\})$ and the three (ENP)s for $\R^2$, obtained from it, given by   $(\{w_j\},\{1/j\}), \,(\{x_j\},\{1/j\}),$ and $(\{y_j\},\{1/j\}),$ where $w_j=(j!,1),$ $ x_1=(1,1)$ and $x_j=(j!,3)$ for $j\geq 2,$ and $y_j=(j!,1/j).$  We see that  the finitely generated groups $G_w$ and  $G_x$ are not topologically isomorphic to $G_y$, since it is not finitely generated.   Of course $G_w$ and $G_x$ are both $\Z \times \Z$ as groups, but  the topologies are distinct.  If not, then both of the sequences $\{(j!,3)\}$ and $\{(j!,1)\}$ would converge to $(0,0)$ in each group, and thus their difference, $\{(0,2)\},$ would converge to $(0,0).$  Since the spaces are Hausdorff, this cannot happen.  Thus on a global scale, there is considerable diversity among the groups formed by (ENP)s.  The remarkable fact is that there is not any topological diversity on a local scale.

Before proceeding with our main results, we need a lemma.  Although the proof of part (iii) of the next lemma is not difficult, we feel it is prudent to write out all of the details, since ``projection'' is not generally continuous with topologies generated by (SNP)s or (ENP)s \cite[pg. 52]{lclisom}.  The reason that it is continuous below is because of a special choice of sequences.

\begin{lem}\label{lem2}
 Let $(\{v_j\},\{p_j\})$ be a (SNP) for $\R^n,\; n>1,$ and let $\{b_j\}=\{(0,0,\ldots,0)\}$ be a sequence in $\R^{n-1}.$
\begin{enumerate}
\renewcommand{\labelenumi}{(\roman{enumi})}
\item If $\{x_j\}$ is defined as in \ref{lem}, then $(\{(x_j,b_j)\},\{p_j\})$ is a (SNP) for  $\R^n.$
\item If $(\{v_j\},\{p_j\})$ satisfies the hypotheses of \ref{thm:isom}, then $(\{(x_j,b_j)\},\{p_j\})$ does as well.
\item The group generated by the pair $(\{x_j\},\{p_j\})$ is topologically isomorphic to the group generated by the pair $(\{(x_j,b_j)\},\{p_j\})$ and the completions of the respective groups are also topologically isomorphic.

\end{enumerate}
\end{lem}
\begin{proof}
Parts (i) and (ii) follow easily from \ref{lem} and definitions.   

Since $\{b_j\}=\{(0,0,\ldots,0)\}$, we will denote the group generated by $\{(x_j,b_j)\}$ as $G_x\times \{0\}^{n-1}$.  Also let $\sigma_1:G_x\times \{0\}^{n-1}\rightarrow \R$ be the norm generated by the pair   $(\{(x_j,b_j)\},\{p_j\})$ and $\sigma_2:G_x\rightarrow \R$ be the norm generated by the pair $(\{x_j\},\{p_j\})$.  Define $p:(G_x\times \{0\}^{n-1},\sigma_1)\rightarrow (G_x,\sigma_2)$ to be projection.    Suppose that $\{(a_j,0,0,\ldots,0)\}$ converges to the identity in $(G_x\times \{0\}^{n-1},\sigma_1).$  Thus for any $\epsilon>0$ there is a $N\in\N$ such that if $j>N$ then 
$$\sigma_1(a_j,0,0,\ldots,0)=\inf \left\{ \sum{|c_k|p_k}:(a_j,0,0,\ldots,0)=\sum{c_k (x_k,b_k)}\right\}<\epsilon.$$
Thus there are summands so that $(a_j,0,0,\ldots,0)=\sum{c_k (x_k,b_k)}$ and $\sum{|c_k|p_k}<\epsilon.$  This is true if and only if there are summands such that $\sum{|c_k|p_k}<\epsilon$  and $a_j=\sum{c_k x_k}$ if and only if 
$$\sigma_2(a_j)=\inf \left\{ \sum{|c_k|p_k}:a_j=\sum{c_k x_k}\right\}\\ <\epsilon$$
 if and only if $\{a_j\}$ converges to the identity in $(G_x,\sigma_2)$.  Hence $(G_x\times \{0\}^{n-1},\sigma_1)$ is topologically isomorphic to $(G_x,\sigma_2)$.  Since $\sigma_1$ and $\sigma_2$ induce invariant metrics (and thus generate the respective uniformities) we can apply \cite[Proposition 5, pg 246]{bb} and conclude that the completions are also topologically isomorphic.
\end{proof}

\begin{thm}\label{thm:rn}
Let $(\{v_j\},\{p_j\})$ be a (SNP) for $\R^n$ that satisfies the hypotheses of \ref{thm:isom}.  If $(G_v,\sigma_v)$ is the topological subgroup generated by $\{v_j\}$, then $\mathcal{C}(G_v,\sigma_v)$ is one dimensional, locally totally disconnected, and not locally compact.

\end{thm}
\begin{proof}
By the previous lemma, $(\{(x_j,0,0,\ldots,0)\},\{p_j\})$ is a (SNP) for  $\R^n$ satisfying the hypotheses of \ref{thm:isom}.  Thus $\mathcal{C}(G_v,\sigma_v)$ is locally isometric to the completion of the group generated by the pair $(\{(x_j,0,0,\ldots,0)\},\{p_j\})$, which, again by the previous lemma, is topologically isomorphic to the completion of the group generated by the pair  $(\{x_j\},\{p_j\})$.  Since $(\{x_j\},\{p_j\})$ is a (NNP) and  $\sum_{j=1}^{\infty}p_j\frac{x_{j+1}}{x_j}$ is unbounded, we can apply \ref{thm:nnp} to conclude that  $\mathcal{C}(\Z,\sigma_x)$ is not locally compact, is one dimensional, and is totally disconnected. We have the statement.
\end{proof}

Our final results show that the groups generated by (ENP)s have the same local structure as the groups generated by (SNP)s.

\begin{lem}\label{lem4}
Let $(\{v_j\},\{p_j\})$ be a (SNP) for $\R^n$ that satisfies the hypotheses of \ref{thm:isom}, and $(\{(v_j,q_j)\},\{p_j\})$ be a corresponding (ENP) for $\R^{n+m}.$ 
 The  group generated by the pair $(\{(v_j,q_j)\},\{p_j\})$ is locally homeomorphic to the group generated by the pair $(\{v_j\},\{p_j\})$, and their respective completions are locally homeomorphic.
\end{lem}
\begin{proof}
Consider the (ENP) for $\R^{n+m}$ given by $(\{(v_j,0,0,\ldots,0)\},\{p_j\}).$  It satisfies the hypotheses of \ref{thm:isom} since $(\{v_j\},\{p_j\})$ does.  Thus the groups generated by the (ENP)s $(\{(v_j,0,0,\ldots,0)\},\{p_j\})$  and $(\{(v_j,q_j)\},\{p_j\})$ are locally isometric.  As in the proof of the previous lemma, the groups generated by $(\{(v_j,0,0,\ldots,0)\},\{p_j\})$ and $(\{v_j\},\{p_j\})$ are topologically isomorphic.  Therefore the groups generated by $(\{v_j\},\{p_j\})$ and $(\{(v_j,q_j)\},\{p_j\})$ are locally homeomorphic.
\end{proof}

\begin{thm}
Let $(\{v_j\},\{p_j\})$ be a (SNP) for $\R^n$ that satisfies the hypotheses of \ref{thm:isom}, and $(\{(v_j,q_j)\},\{p_j\})$ be a corresponding (ENP) for $\R^{n+m}.$ 
If $(G_v,\sigma_v)$ is the topological subgroup generated by $\{(v_j,q_j)\}$, then $\mathcal{C}(G_v,\sigma_v)$ is one dimensional, locally totally disconnected, and not locally compact.

\end{thm}
\begin{proof}
Apply the \ref{lem4} to the group generated by the (ENP).  We then know that this group is locally homeomorphic to a group generated by a (SNP).  Thus \ref{thm:rn} gives the result.
\end{proof}

\section{Concluding Remarks}
As we have seen, this paper deals with the local topological structure of certain groups, and we have exploited the author's and Nienhuys's analysis to understand certain aspects of  this local structure.  Understanding the global structure of  these groups is proving particularly difficult.  We are currently considering certain cases.    In particular, we are considering groups where the main condition of \ref{thm:isom} is an equality (at least asymptotically), i.e.,  we consider $\lim_{j\rightarrow \infty}{p_j\left\lceil\frac{\|v_{j+1}\|}{\|v_j\|}-1\right\rceil}=1$ for example.  Even though it does not simplify the situation much, it does allow us to have more ``control'' over the sequences and groups they generate. 

In \cite{lclisom} we were primarily concerned with the structure of the subgroups and quotient groups of $\R^n$ obtained using the (SNP) $(\{v_j\},\{p_j\})$ and the norm $\psi:\R^n\rightarrow \R$ given by $\psi(x)=\inf\left\{\sum{|c_j|p_j}+\|x-\sum{c_j v_j}\|\right\}$.  Although $\psi$ certainly resembles $\sigma$, it is different in a couple of important ways.  First, it is defined on the full group $\R^n$ whereas $\sigma$ is defined on a proper subgroup of $\R^n$. Second, it induces the same topology on $G_v$ as $\sigma$, if $\{v_j\}$ generates (in the standard sense) a discrete subgroup, but it gives a weaker group topology on $G_v$, if $\{v_j\}$ generates a nondiscrete subgroup of $\R^n$ \cite[Proposition 12]{lclisom}.   In either case, $\psi$ has the added attraction of always giving group topologies that are weaker than the standard topologies.  (Thus this norm is related to the study of immersions of Lie groups.)  
We plan to address the structure of these groups in a future paper.

\end{document}